\documentclass{amsart}
\usepackage{latexsym,amsxtra,amscd,ifthen}
\usepackage{amsfonts}
\usepackage{pifont}
\usepackage{verbatim}
\usepackage{amsmath}
\usepackage{graphicx}
\usepackage{amsthm}
\usepackage{amssymb}
\usepackage{url,color}
\usepackage[arrow,matrix]{xy}
\usepackage{pict2e,epic,eepic,pstricks}
\usepackage{multirow}
\usepackage[colorlinks=true,final,hyperindex=true,linkcolor=black,citecolor=black,urlcolor=black,filecolor=black]{hyperref}
\usepackage{bbm}
\usepackage{float}
\usepackage{mathtools}
\usepackage{cancel}
\usepackage{anysize}\marginsize{25mm}{25mm}{30mm}{38mm}
\allowdisplaybreaks[4]
\newtheorem{thm}{Theorem}[section]

\newtheorem{lem}[thm]{Lemma}

\newtheorem{prop}[thm]{Proposition}
\newtheorem{prop-def}[thm]{Proposition-Definition}
\theoremstyle{definition}
\newtheorem{Def}[thm]{Definition}
\theoremstyle{remark}
\newtheorem{rmk}[thm]{\bf Remark}
\newtheorem{exm}[thm]{\bf Example}
\numberwithin{equation}{section}

\newcommand{\lr}[1]{\langle #1 \rangle}

\def \Ga{\Gamma}
\def \ga{\gamma}
\def \la{\lambda}
\def \si{\sigma}
\def \vp{\varphi}
\def \k{\Bbbk}
\def \H{\mathbb{H}}
\def \M{\mathbb{M}}

\def \P{\mathbb{P}}
\def \S{\mathbb{S}}
\def \Z{\mathbb{Z}}
\def \K{\mathcal{K}}
\def \gd{\operatorname{Gdeg}}
\def \Id{\operatorname{Id}}
\def \sp{\operatorname{span}}
\def \T{\operatorname{T}}
\def \id{\mathrm{id}}
\def \es{\emptyset}

\begin{document}
	
\title{Superalgebras and Algebras with Involution: Classifying Cubic Codimension Sequences}
	
\author{Yan-Hong Bao$^*$}
\address{(Bao) School of Mathematical Sciences, Anhui University, Hefei 230601, China}
\email{baoyh@ahu.edu.cn}
\thanks{$^*$Corresponding author: Yan-Hong Bao,\href{mailto:baoyh@ahu.edu.cn}{baoyh@ahu.edu.cn}}
	
\author{Jiang-Nan Xu}
\address{(Xu) School of Finance and Mathematics, Huainan Normal University, Huainan 232038, China}
\email{jnxu1@outlook.com}
	
\author{Yuan-Feng Zhang}
\address{(Zhang)  School of Mathematical Sciences, Anhui University, Hefei 230601, China}
\email{yfzhang@stu.ahu.edu.cn}

	

\begin{abstract} 
	A $\vp$-algebra is either a superalgebra or an algebra with involution. In this paper, we study $\T^\vp$-ideals associated with unital $\vp$-algebras whose $\vp$-codimension sequence exhibits cubic polynomial growth. As a consequence, we obtain a complete classification of all $\vp$-codimension sequences of cubic growth for unital $\vp$-algebras. Furthermore, we explicitly determine a minimal-degree multilinear generator for every $\T^\vp$-ideal associated with unital $\vp$-algebras whose $\vp$-codimension growth is at most quadratic.
\end{abstract}
	
\maketitle
\setcounter{section}{-1}

\section{Introduction} \label{sec0} 
	Throughout this paper, let $\k$ be a field of characteristic zero. In the early 1950s, Specht \cite{Sp}  conjectured that every proper T-ideal of free algebra is finitely generated as a T-ideal. Since then, classifying T-ideals and determining their explicit generators for specific classes of PI-algebras have remained central problems in PI theory. Although Kemer \cite{Ke87} confirmed Specht's conjecture for fields of characteristic zero in 1987, explicitly describing the full structure of a $\T$-ideal remains a highly difficult problem in general.
	
	To overcome this difficulty, one may study the asymptotic behavior of T-ideals via numerical growth invariants. A fundamental tool in this regard is the codimension sequence, introduced by Regev \cite{Re} in 1972, who also proved that this sequence is exponentially bounded. Kemer \cite{Ke79} showed that the codimension sequence of any PI-algebra can only grow polynomially or exponentially, with no intermediate behavior. The classification of $\T$-ideals with low polynomial codimension growth has also been studied in \cite{BFXYZZZ,BYZ,GM,GMP,GMZ,OV}.

	Let $A$ be an algebra over the field $\k$, and let $\vp$ be an automorphism or antiautomorphism of $A$ satisfying $\vp^2 = \id_A$. Recall that a \emph{$\vp$-algebra} is a pair $(A,\vp)$ consisting of an algebra $A$ equipped with such a map $\vp$. When $\vp$ is an automorphism, its $+1$ and $-1$ eigenspaces define a $\Z_2$-grading on $A$, so that $(A,\vp)$ can be viewed as a \emph{superalgebra}, and when $\vp$ is an antiautomorphism, it is called an involution and $(A,\vp)$ is termed a \emph{$*$-algebra}. The pair $(A,\vp)$ is often abbreviated simply as $A$ when $\vp$ is clear from context.

	Similar to the codimension sequence of ordinary PI-algebras, one can define the $\vp$-codimension sequence for $\vp$-algebras. Specifically, given a $\vp$-algebra $A$, denote by $\Id^\vp(A)$ its $\T^\vp$-ideal of $\vp$-identities, and by $V_n^\vp$ the space of multilinear polynomials of degree $n$ in the variables $y_1,z_1,\dots,y_n,z_n$. The \emph{$\vp$-codimension sequence} $\{c_n^\vp(A)\}_{n \ge 1}$ of $A$ is defined by
	\[
	c_n^\vp(A) = \dim_\k \frac{V_n^\vp}{V_n^\vp \cap \Id^\vp(A)}.
	\]
	In the automorphism case this is called the \emph{graded codimension sequence}, and in the antiautomorphism case the \emph{$*$-codimension sequence}. The sequence is said to have polynomial growth $n^k$ if there exists a constant $q > 0$ such that
	\[
	c_n^\vp(A) = q n^k + O(n^{k-1}), \quad  \text{as} \quad  n \to \infty.
	\]
	
	The classification of $\vp$-algebras via $\vp$-codimension growth was studied in the linear case by Giambruno, La Mattina and Misso \cite{GMM}, who treated superalgebras and determined the corresponding $\T^\vp$-ideals, and the $*$-algebra case was settled analogously in \cite{MM}. Minimal unital $\vp$-varieties of polynomial $\vp$-codimension growth, together with their $\T^\vp$-ideals, were characterized by Gouveia, dos Santos and Vieira \cite{GSV} in 2020. Building on these results, Bessades, dos Santos, Santos and Vieira \cite{BSSV} completely classified unital $\vp$-algebras with quadratic $\vp$-codimension growth in 2021, explicitly determining both the algebras and their associated $\T^\vp$-ideals: they exhibited $16$ superalgebras giving rise to $6$ distinct graded codimension sequences, and $15$ $*$-algebras corresponding to $5$ distinct $*$-codimension sequences.
	
	Just as proper polynomials are crucial for ordinary unital PI-algebras, so are $Y$-proper polynomials for $\vp$-algebras. Denote by $\Ga_n^\vp$ the space of multilinear $Y$-proper $\vp$-polynomials of degree $n$. By Theorem 4.3.3 of \cite{Dr}, the $\T^\vp$-ideal $\Id^\vp(A)$ is completely determined by its intersections with the spaces $\Ga_n^\vp$. Thus, the study of $\Id^\vp(A)$ reduces to the sequence
	\[
	(\Ga_n^\vp \cap \Id^\vp(A))_{n \geq 1} = (\Ga_1^\vp \cap \Id^\vp(A), \Ga_2^\vp \cap \Id^\vp(A), \dots).
	\]
	Let $\Ga_0^\vp=\sp_\k\{1\}$. Recall that the \emph{proper $\vp$-codimension sequence} $\{\ga_n^\vp(A)\}_{n\geq 0}$ is defined by
	\[
	\ga_n^\vp(A) \coloneqq \dim_\k \frac{\Ga_n^\vp}{\Ga_n^\vp \cap \Id^\vp(A)}.
	\]
	By Theorem 2.3 of \cite{DG}, the $\vp$-codimension sequence and the proper $\vp$-codimension sequence satisfy the binomial relation
	\begin{equation} \label{0.0.1} \tag{E0.0.1}
		c_n^\vp (A) = \sum_{i=0}^{n} \binom{n}{i} \ga_i^\vp(A).
	\end{equation}
	Hence, classifying all $\T^\vp$-ideals with polynomial growth $n^k$ for unital $\vp$-algebras further reduces to considering the sequence 
	\[
	(\Ga_1^\vp \cap \Id^\vp(A), \dots, \Ga_k^\vp \cap \Id^\vp(A),\Ga_{k+1}^\vp, \Ga_{k+2}^\vp, \dots).
	\]
	In \cite{BFXYZZZ}, the authors introduced \emph{$k$-admissible sequences} to characterize such truncated sequences of the first $k$ terms. A similar definition is given in the present paper (see Definition~\ref{def2.1}).
	
	In this paper, we study $\T^\vp$-ideals with cubic $\vp$-codimension growth for unital $\vp$-algebras. By classifying the associated $3$-admissible sequences, we obtain a partial classification of these ideals in the superalgebra setting, yielding $59$ distinct graded codimension sequences (see Theorem~\ref{thm2.10}). For $*$-algebras, we obtain a complete classification of such ideals, resulting in $54$ distinct $*$-codimension sequences (see Theorem~\ref{thm2.13}).
	
	In \cite{BXYZZ}, using the finite basis property of T-ideals, the authors proved that every T-ideal associated with a unital algebra is generated by a single element. Subsequently, in \cite{BFXYZZZ}, the authors studied a minimal-degree multilinear generator for the T-ideal associated with a unital PI-algebra, called the \emph{generating identity} of the algebra, and determined the generating identities for unital PI-algebras whose codimension growth is at most quartic. We consider a similar problem for unital $\vp$-algebras, and determine the generating identities for unital $\vp$-algebras whose $\vp$-codimension growth is at most quadratic (see Proposition~\ref{prop3.6}).
	
	The paper is organized as follows. Section~1 recalls the basic notions and key tools from the PI theory of $\vp$-algebras. Section~2 introduces $k$-admissible sequences and their relation to $\T^\vp$-ideals with $\vp$-codimension growth $n^k$, then classifies all $3$-admissible sequences, thereby yielding the complete classification of cubic $\vp$-codimension sequences. Section~3 proves the existence of generating identities for unital $\vp$-algebras, gives a general method for computing them, and presents the generating identities for $\T^\vp$-ideals with at most quadratic $\vp$-codimension growth.
	
\section{Preliminaries} \label{sec1}
	Throughout this paper, all algebras are assumed to be associative and unital. For a positive integer $n$, let $[n] = \{1,2,\dots,n\}$, and denote by $\S_n$ the symmetric group on $n$ letters.
	
	A \emph{$\vp$-algebra} is a pair $(A, \vp)$ where $\vp$ is an automorphism or antiautomorphism of $A$ satisfying $\vp^2 = \id_A$. That is, for all $a,b \in A$, either $\vp(ab) = \vp(a)\vp(b)$ or $\vp(ab) = \vp(b)\vp(a)$. Let $(A,\vp_A)$ and $(B,\vp_B)$ be $\vp$-algebras. An algebra homomorphism $\psi \colon A \to B$ is said to be a \emph{$\vp$-algebra homomorphism} if $\psi \circ \vp_A = \vp_B \circ \psi$. Let $X = \{x_1, x_2, \dots\}$ be a countable set, and let $\k \lr{X, \vp} = \k \lr{x_1, \vp(x_1), x_2, \vp(x_2), \dots}$ denote the \emph{free $\vp$-algebra}. Recall that the free $\vp$-algebra on $X$ is characterized by the following universal property: given a $\vp$-algebra $(A, \vp_A)$, any map $f \colon X \to A$ can be uniquely extended to a $\vp$-algebra homomorphism $\hat{f} \colon \k \langle X, \vp \rangle \to A$. For each $i \geq 1$, set $y_i = x_i + \vp(x_i)$ and $z_i = x_i - \vp(x_i)$. Then $\k \lr{X, \vp} = \k \lr{Y, Z}$, where $Y = \{y_1, y_2, \dots\}$ and $Z = \{z_1, z_2, \dots\}$. By construction, $\vp(y_i) = y_i$ and $\vp(z_i) = -z_i$. Elements of $Y$ are called symmetric variables, and elements of $Z$ are called skew variables.
	
	Since the base field is of characteristic zero, every $\vp$-algebra $A$ decomposes as a direct sum of subspaces:
	\[
	A = A_\vp^+ \oplus A_\vp^-, \quad \text{where} \quad A_\vp^+ = \{ a \in A \mid \vp(a) = a \}, \quad A_\vp^- = \{ a \in A \mid \vp(a) = -a \}.
	\]
	A polynomial $f(y_1, \dots, y_n, z_1, \dots, z_m) \in \k \lr{Y, Z}$ is called a \emph{$\vp$-identity} of $A$ if 
	\[
	f(a_1, \dots, a_n, b_1, \dots, b_m) = 0
	\] 
	for all $a_1, \dots, a_n \in A_\vp^+$ and $b_1, \dots, b_m \in A_\vp^-$. Denote by $\Id^\vp(A)$ the set of all $\vp$-identities of $A$. Then $\Id^\vp(A)$ is a \emph{$\T^\vp$-ideal} of $\k \lr{Y, Z}$, i.e., an ideal invariant under all $\vp$-algebra endomorphisms of $\k \lr{Y, Z}$.
	
	Let $V_n^\vp$ be the space of multilinear polynomials of degree $n$ in the variables $y_1, z_1, \dots, y_n, z_n$, namely
	\[
	V_n^\vp = \sp_\k \{ v_{\si(1)} \cdots v_{\si(n)} \mid \si \in \S_n, v_i \in \{y_i, z_i\} \text{ for } i \in [n] \}.
	\]
	For a $\vp$-algebra $A$, let
	\[
	V_n^\vp(A) = \frac{V_n^\vp}{V_n^\vp \cap \Id^\vp(A)},
	\qquad
	c_n^\vp(A) = \dim_\k V_n^\vp(A).
	\]
	The sequence $\{c_n^\vp(A)\}_{n \geq 1}$ is called the \emph{$\vp$-codimension sequence} of $A$. It is one of the fundamental invariants describing the growth of the associated $\T^\vp$-ideal.  
	
	Let $\Z_2 = \{1, -1\}$ be the cyclic group of order $2$. The wreath product $\Z_2 \wr \S_n$, known as the hyperoctahedral group of degree $n$ and denoted by $\H_n$, is 
	\[
	\H_n = \{(a_1, \dots, a_n; \si) \mid a_i \in \Z_2, \si \in \S_n\},
	\]
	with multiplication given by $(a_1, \dots, a_n; \si)(b_1, \dots, b_n; \tau) = (a_1 b_{\si^{-1}(1)}, \dots, a_n b_{\si^{-1}(n)}; \si \tau)$. The space $V_n^\vp$ admits a natural left $\H_n$-action defined on symmetric and skew variables by
	\[
	h \cdot y_i = y_{\si(i)} \quad \text{and} \quad h \cdot z_i = a_{\si(i)} z_{\si(i)},
	\]
	for $h = (a_1, \dots, a_n; \si) \in \H_n$. For integers $0 \leq r \leq n$, let $\la \vdash r$, $\mu \vdash (n-r)$ be partitions. We call the pair $(\la, \mu)$ a \emph{bipartition} of $n$, denoted by $(\la, \mu) \vdash (r, n-r)$. The isomorphism classes of irreducible $\k\H_n$-modules are in bijective correspondence with the bipartitions of $n$ \cite[Subsection 10.6]{GZ}. We write $\chi_{\la,\mu}$ for the irreducible $\H_n$-character associated with the bipartition $(\la,\mu)$, and for a left $\k \H_n$-module $M$, we denote its $\H_n$-character by $\chi_n(M)$, or simply by $\chi(M)$ when $n$ is clear from context. For a $\vp$-algebra $A$, the $\H_n$-character $\chi_n(V_n^\vp(A))$ is called the \emph{$n$-th $\vp$-cocharacter} of $A$. By complete reducibility, we have the decomposition
	\begin{equation} \label{1.0.1} \tag{E1.0.1}
		\chi_n(V_n^\vp(A)) = \sum_{(\la,\mu) \vdash n} m_{\la,\mu} \, \chi_{\la,\mu} = \sum_{r=0}^{n} \, \sum_{(\la,\mu) \vdash (r,n-r)}^{} m_{\la,\mu} \, \chi_{\la,\mu},
	\end{equation}
	where $m_{\lambda,\mu}$ denotes the multiplicity of the irreducible character $\chi_{\lambda,\mu}$.
	
	Let $V_{r,n-r}^\vp$ be the subspace of $V_n^\vp$ spanned by all multilinear monomials with the first $r$ symmetric variables $y_1, \dots, y_r$ and the remaining $n-r$ skew variables $z_{r+1}, \dots, z_n$, namely,
	\[
	V_{r,n-r}^\vp = \sp_\k \{ v_{\si(1)} \cdots v_{\si(n)} \mid \si \in \S_n, \, v_i = y_i \ (1 \leq i \leq r), \, v_j = z_j \ (r+1 \leq j \leq n) \}.
	\]
	The space $V_{r,n-r}^\vp$ carries a natural left $\S_r \times \S_{n-r}$-action. The irreducible characters of $\S_r \times \S_{n-r}$ are exactly the tensor products of irreducible $\S_r$-characters and irreducible $\S_{n-r}$-characters. The isomorphism classes of irreducible $\k[\S_r \times \S_{n-r}]$-modules correspond bijectively to bipartitions $(\lambda,\mu)$ of $n$. We write $\chi_\la \otimes \chi_\mu$ for the irreducible $\S_r \times \S_{n-r}$-character associated with the bipartition $(\la,\mu)$, and for a left $\k[ \S_r \times \S_{n-r}]$-module $M$, we denote its $\S_r \times \S_{n-r}$-character by $\chi_{r,n-r}(M)$. For a $\vp$-algebra $A$, let
	\[
	V_{r,n-r}^\vp(A) = \frac{V_{r,n-r}^\vp}{V_{r,n-r}^\vp \cap \Id^\vp(A)},
	\qquad
	c_{r,n-r}^\vp(A) = \dim_\k V_{r,n-r}^\vp(A).
	\]
	By complete reducibility, for each $0 \leq r \leq n$ we have the decomposition
	\begin{equation} \label{1.0.2} \tag{E1.0.2}
		\chi_{r,n-r}(V_{r,n-r}^\vp(A)) = \sum_{(\la,\mu) \vdash (r,n-r)} \tilde{m}_{\la,\mu} \, \chi_\la \otimes \chi_\mu,
	\end{equation}
	where $\tilde{m}_{\la,\mu}$ denotes the multiplicity of the irreducible character $\chi_\la \otimes \chi_\mu$.
	
	\begin{thm}\cite[Theorem 1.3]{DG}\label{thm1.1}
		The multiplicities $m_{\la,\mu}$ and $\tilde{m}_{\la,\mu}$ appearing in equations \eqref{1.0.1} and \eqref{1.0.2} are equal. Moreover, we have the following codimension relation:
		\[
		c_n^\vp(A) = \sum_{r=0}^{n} \binom{n}{r} c_{r,n-r}^\vp(A).
		\]
	\end{thm}
	Hence, the $\k\H_n$-module structure of $V_n^\vp \cap \Id^\vp(A)$ is completely determined by the $\k[\S_r \times \S_{n-r}]$-module structures of the homogeneous components $V_{r,n-r}^\vp \cap \Id^\vp(A)$ for all $0 \leq r \leq n$.
	
	For unital $\vp$-algebras, the study of $\vp$-identities can be reduced to that of $Y$-proper $\vp$-identities. Let $B(Y)$ denote the subalgebra of the free $\vp$-algebra $\k \lr{Y, Z}$ generated by the variables in $Z$ and all commutators of elements in $Y \cup Z$. Elements of $B(Y)$ are called \emph{$Y$-proper $\vp$-polynomials} and admit a normal form: 
	\[
	z_{i_1} \cdots z_{i_l} [v_{j_1},\dots,v_{j_m}] \cdots [v_{k_1},\dots,v_{k_n}],
	\]
	where $i_1 < \dots < i_l$, $v_i \in Y \cup Z$, and the iterated commutator is defined recursively by $[v_1,\dots,v_{t-1},v_t] = [[v_1,\dots,v_{t-1}],v_t]$. For each integer $n \geq 0$, set $\Ga_n^\vp = B(Y) \cap V_n^\vp$, with $\Ga_0^\vp = \sp_\k \{1\}$. Lemma 2.1 of \cite{MMM} computes the dimension of $\Ga_n^\vp$ as
	\[
	\dim_{\k} \Ga_n^\vp = n! \sum_{i=0}^{n} 2^{n-i} \frac{(-1)^i}{i!}.
	\]
	By \cite[Theorem 4.3.3]{Dr}, the $\T^\vp$-ideal $\Id^\vp(A)$ is completely determined by its $Y$-proper multilinear components. Consequently, the study of $\Id^\vp(A)$ reduces to the sequence $(\Ga_n^\vp \cap \Id^\vp(A))_{n \geq 1}$.
	
	The \emph{proper $\vp$-codimension sequence} $\{\ga_n^\vp(A)\}_{n\geq 0}$ is given by
	\[
	\ga_n^\vp(A) \coloneqq \dim_\k \frac{\Ga_n^\vp}{\Ga_n^\vp \cap \Id^\vp(A)}.
	\]
	Combining the binomial relation \eqref{0.0.1} and Remark 2.2 in \cite{BSSV}, we obtain an equivalent characterization of polynomial growth $n^k$ for the \(\varphi\)-codimension sequence.
	\begin{lem} \label{lem1.2}
		The $\vp$-codimension sequence $c_n^\vp(A)$ has polynomial growth $n^k$ if and only if $\Ga_k^\vp \not\subset \Id^\vp(A)$ and $\Ga_i^\vp \subset \Id^\vp(A)$ for all $i \geq k+1$.
	\end{lem}
	\begin{proof}
		By \eqref{0.0.1} and Remark~2.2 of \cite{BSSV}, $\ga_i^\vp(A)\ge 0$ and $\ga_i^\vp(A)=0$ if and only if $\Ga_i^\vp\subset\Id^\vp(A)$. If $\Ga_k^\vp\not\subset\Id^\vp(A)$ and $\Ga_i^\vp\subset\Id^\vp(A)$ for all $i>k$, then $\ga_k^\vp(A)>0$ and
		\[
		c_n^\vp(A)=\sum_{i=0}^k \ga_i^\vp(A)\binom{n}{i}
		=\frac{\ga_k^\vp(A)}{k!}n^k+O(n^{k-1}).
		\]
		Conversely, assume $c_n^\vp(A)=q n^k+O(n^{k-1})$ with $q>0$. If $\ga_i^\vp(A)>0$ for some $i>k$, then
		\[
		c_n^\vp(A)\ge \ga_i^\vp(A)\binom{n}{i}
		=\frac{\ga_i^\vp(A)}{i!}n^i+O(n^{i-1}),
		\]
		which grows faster than $n^k$, a contradiction. Hence $\ga_i^\vp(A)=0$ for all $i>k$. Since the leading coefficient is $\ga_k^\vp(A)/k!=q>0$, it follows that $\ga_k^\vp(A)>0$. The two containments follow.
	\end{proof}
	The subspace $\Ga_n^\vp \subset V_n^\vp$ is invariant under the action of $\H_n$, and thus it is a left $\k\H_n$-module. In \cite[Section 4]{GSV}, the authors use representation theory of the general linear group to decompose the $\H_n$-character of $\Ga_n^\vp$ into irreducible characters for small $n$:
	\begin{equation*}
		\begin{split}
			\chi(\Ga_1^\vp) &= \chi_{\es, (1)}, \\
			\chi(\Ga_2^\vp) &= \chi_{(1^2),\es} + \chi_{(1),(1)} + \chi_{\es,(2)} + \chi_{\es, (1^2)}, \\
			\chi(\Ga_3^\vp) &= \chi_{(2,1),\es} + \chi_{(2),(1)} + 2\chi_{(1^2),(1)} + 2\chi_{(1),(2)} + 2\chi_{(1),(1^2)} + \chi_{\es, (3)} + 2\chi_{\es, (2,1)} + \chi_{\es, (1^3)}.
		\end{split}
	\end{equation*}
	
	In the classical setting of ordinary unital associative algebras, a Specht basis for the proper polynomial space $\Ga_n$ is constructed in \cite[Theorem 4.3.9]{Dr}. An analogous result holds for unital $\vp$-algebras. For degree $1$, the Specht basis of $\Ga_1^\vp$ is given by $\{z_1\}$. For degree $2$, the Specht basis of $\Ga_2^\vp$ consists of the elements
	\[
	[y_2,y_1],\ [z_2,y_1],\ [y_2,z_1],\ [z_2,z_1],\ z_1z_2,
	\]
	and $\Ga_2^\vp$ admits the following decomposition as a left $\k\H_2$-module:
	\[
	\Ga_2^\vp = V_{(1^2), \es} \oplus V_{(1), (1)} \oplus V_{\es, (2)} \oplus V_{\es, (1^2)},
	\]
	where
	\[
		V_{(1^2), \es} = \sp_\k \{ [y_2,y_1] \}, 
		V_{(1), (1)} = \sp_\k \{ [z_2,y_1], [y_2,z_1] \}, 
		V_{\es, (2)} = \sp_\k \{ z_1\circ z_2 \}, 
		V_{\es, (1^2)} = \sp_\k \{ [z_2,z_1] \},
	\]
	and $z_1 \circ z_2 = z_1 z_2 + z_2 z_1$ denotes the Jordan product. For degree $3$, a Specht basis of $\Ga_3^\vp$ is listed below:
	\[
	\begin{array}{rrrrrr}
		\relax [y_3,y_1,y_2] & [y_3,y_2,y_1] & & & & \\
		\relax [z_3,y_1,y_2] & [z_3,y_2,y_1] & [y_3,y_1,z_2] & [y_3,z_2,y_1]  & [y_3,z_1,y_2] & [y_3,y_2,z_1] \\
		z_3[y_2,y_1]         & z_2[y_3,y_1]  & z_1[y_3,y_2]  & & &\\
		\relax [z_3,y_1,z_2] & [z_3,z_2,y_1] & [z_3,z_1,y_2] & [z_3,y_2,z_1] & [y_3,z_1,z_2] & [y_3,z_2,z_1]\\ 
		z_2[z_3,y_1]         & z_3[z_2,y_1]  & z_1[z_3,y_2]  & z_3[y_2,z_1]  & z_1[y_3,z_2]  & z_2[y_3,z_1] \\
		z_1z_2z_3            & z_1[z_3,z_2]  & z_2[z_3,z_1]  & z_3[z_2,z_1]  & [z_3,z_1,z_2] & [z_3,z_2,z_1]
	\end{array}
	\]
	The irreducible $\k\H_n$-submodules of $\Ga_n^\vp$ can be obtained by inducing irreducible $\k[\S_r \times \S_{n-r}]$-submodules of the homogeneous components $\Ga_n^\vp \cap V_{r,n-r}^\vp$. For example, let $k_1,k_2 \in \k$ be scalars, not both zero, and let $f_1 = k_1([z_3,y_1,y_2] - [z_3,y_2,y_1]) + k_2 z_3[y_2,y_1]$. Note that $V_{(1^2),(1)} = \sp_\k\{ f_1 \}$ is an irreducible $\S_2 \times \S_1$-submodule of $\Ga_3^\vp$ corresponding to the irreducible character $\chi_{(1^2)} \otimes \chi_{(1)}$. Let
	\[
	\begin{aligned}
		f_2 &= k_1([z_2,y_1,y_3] - [z_2,y_3,y_1]) + k_2 z_2[y_3,y_1]= k_1[y_3,y_1,z_2]+k_2 z_2[y_3,y_1], \\
		f_3 &= k_1([z_1,y_3,y_2] - [z_1,y_2,y_3]) + k_2 z_1[y_2,y_3]= -k_1[y_3,y_2,z_1]-k_2 z_1[y_3,y_2].
	\end{aligned}
	\]
	Then the $\k\H_3$-module generated by $V_{(1^2),(1)}$ is
	\[
	\k\H_3 \cdot V_{(1^2),(1)} =\sp_\k \{ f_1, f_2, f_3 \},
	\]
	which is an irreducible $\k\H_3$-submodule of $\Ga_3^\vp$ corresponding to the irreducible character $\chi_{(1^2),(1)}$.

\section{Cubic $\vp$-codimension sequences of unital $\vp$-algebras} \label{sec2}
	The goal of this section is to classify all $3$-admissible sequences up to componentwise isomorphism of $\k\H_3$-modules, thereby determining the associated $\T^\vp$-ideals and yielding the complete list of cubic $\vp$-codimension sequences.

\subsection{Admissible sequences}
	Let $A$ be a unital $\vp$-algebra, and let $\Id^\vp(A)$ be its associated $\T^\vp$-ideal. By Lemma~\ref{lem1.2}, if the $\vp$-codimension sequence of $A$ has polynomial growth $n^k$, then $\Id^\vp(A)$ is completely determined by the sequence
	\[
	(\Ga_1^\vp \cap \Id^\vp(A), \dots, \Ga_k^\vp \cap \Id^\vp(A), \Ga_{k+1}^\vp, \Ga_{k+2}^\vp, \dots).
	\]
	Hence, classifying all such $\T^\vp$-ideals reduces to determining the first $k$ terms of this sequence.
	
\begin{Def}\label{def2.1}
	Let $s$ and $k$ be integers with $1 \leq s \leq k$. For each $i = s, \dots, k$, let $M_i$ be a left $\k\H_i$-submodule of $\Ga_i^\vp$ satisfying $\vp(M_i) \subset M_i$, and let $\M$ denote the sequence $(M_s, M_{s+1}, \dots, M_k)$. Then $\M$ is called a \emph{$k$-admissible sequence} if the following conditions hold:
	\begin{enumerate}
		\item[(1)] If $s < k$, then $M_s \neq 0$;
		\item[(2)] $M_k$ is a proper submodule of $\Ga_k^\vp$;
		\item[(3)] For every integer $j$ with $s+1 \leq j \leq k$,
		\[
		\bigl( \langle M_s \rangle_{\T^\vp} + \dots + \langle M_{j-1} \rangle_{\T^\vp} \bigr) \cap \Ga_j^\vp \subset M_j,
		\]
		where $\langle M_i \rangle_{\T^\vp}$ denotes the $\T^\vp$-ideal generated by $M_i$.
	\end{enumerate}
	For brevity, we use the following terminology: $\M = (M_k)$ is called a \emph{$k$-admissible singleton}; $\M = (M_{k-1}, M_k)$ is called a \emph{$k$-admissible pair}; and $\M = (M_{k-2}, M_{k-1}, M_k)$ is called a \emph{$k$-admissible triple}. A $\vp$-algebra $A$ is said to be associated with a $k$-admissible sequence $\M=(M_s,\dots,M_k)$ if $\Ga_i^\vp \cap \Id^\vp(A) = 0$ for $1 \le i < s$, $\Ga_i^\vp \cap \Id^\vp(A) = M_i$ for $i=s,\dots,k$, and $\Ga_i^\vp \cap \Id^\vp(A) = \Ga_i^\vp$ for $i>k$. The trivial case $k=0$ corresponds to the empty admissible sequence, meaning that $\Id^\vp(A)$ contains all $Y$-proper polynomials of degree $\ge 1$ and $c_n^\vp(A)=1$ for all $n$.
\end{Def}

\begin{rmk}\label{rmk2.2}
	Every $k$-admissible sequence $\M=(M_s,\dots,M_k)$ determines a unique $\T^\vp$-ideal
	\[
	I_\M=\bigl\langle M_s,\dots,M_k,\ \Ga_{k+1}^\vp,\Ga_{k+2}^\vp,\dots \bigr\rangle_{\T^\vp}.
	\]
	For $j<s$, there are no generators of degree $\le j$, and generators of higher degree cannot produce nonzero $Y$-proper polynomials of degree $j$: every symmetric variable in a $Y$-proper polynomial lies inside some commutator, hence vanishes when replaced by $1$, while skew variables cannot be substituted by $1$. For $s\le j\le k$, condition (3) ensures that contributions from lower-degree generators are contained in $M_j$, and the degree-$j$ generators span exactly $M_j$ as a $\k\H_j$-module. Therefore,
	\[
	I_\M\cap \Ga_j^\vp=
	\begin{cases}
		0, & j<s,\\
		M_j, & s\le j\le k,\\
		\Ga_j^\vp, & j>k.
	\end{cases}
	\]
	Thus $I_\M$ is a proper $\T^\vp$-ideal, and $\k\langle Y,Z\rangle/I_\M$ is a unital $\vp$-algebra associated with $\M$. Conversely, every $\T^\vp$-ideal of growth $n^k$ arises from its first $k$ proper components, so $k$-admissible sequences parametrize such ideals.
\end{rmk}
	
	Expressed in terms of $k$-admissible sequences, relation \eqref{0.0.1} can be restated as follows.
	
\begin{lem}\label{lem2.3}
	Let $A$ be a unital $\vp$-algebra associated with the $k$-admissible sequence $\M=(M_s,\dots,M_k)$. Then
	\[
	c_n^\vp(A) = \sum_{i=0}^{s-1} \dim_\k \Ga_i^\vp \binom{n}{i}+ \sum_{i=s}^{k}(\dim_\k \Ga_i^\vp - \dim_\k M_i) \binom{n}{i}.
	\]
\end{lem}
\begin{proof}
	By \eqref{0.0.1}, the coefficients are the dimensions of the quotient spaces $\Ga_i^\vp/(\Ga_i^\vp \cap \Id^\vp(A))$. For the $k$-admissible sequence $\M$, this quotient has dimension $\dim_\k \Ga_i^\vp$ when $i<s$, $\dim_\k \Ga_i^\vp - \dim_\k M_i$ when $s\le i\le k$, and $0$ when $i>k$. Substituting into \eqref{0.0.1} yields the formula.
\end{proof}

\begin{exm}\label{exm2.4}
	Let $A$ be a unital superalgebra whose graded codimension sequence has quadratic growth. There is a unique $2$-admissible pair of the form
	\[
	\P=(\Ga_1^\vp = \sp_\k \{ z_1 \}, V_{(1),(1)} \oplus V_{\es,(2)} \oplus V_{\es,(1^2)}).
	\]
	To see this, recall that by the substitution property of the  $\T^\vp$-ideal, a skew variable $z$ can be replaced by $yz$ or $zy$, but cannot be replaced by the symmetric element $yy$. By the two-sided ideal property, $z_1$ also produces $z_1z_2$ and $z_2z_1$. From the Specht basis of $\Ga_2^\vp$, these elements generate $[y_2,z_1]$, $[z_2,y_1]$, $z_1 \circ z_2$, and $[z_2,z_1]$, which precisely span $V_{(1),(1)} \oplus V_{\es,(2)} \oplus V_{\es,(1^2)}$.
\end{exm}
	
	The classification of quadratic $\vp$-codimension sequences for unital $\vp$-algebras is given in \cite[Theorem 5.1]{BSSV}. We reformulate this result in the language of $2$-admissible sequences.
	
\begin{prop}\cite[Theorem 5.1]{BSSV}\label{prop2.5}
	Let $A$ be a unital $\vp$-algebra whose $\vp$-codimension sequence $c_n^\vp(A)$ has quadratic growth.
		
	\begin{enumerate}
		\item[(1)] If $\vp$ is an antiautomorphism, then $A$ is associated with a $2$-admissible singleton $\M=(M_2)$, where $M_2$ is a proper $\k \H_2$-submodule of $\Ga_2^\vp = V_{(1^2), \es} \oplus  V_{(1),(1)} \oplus V_{\es,(2)} \oplus V_{\es,(1^2)}$. In this case,
		\[
		c_n^\vp(A)= 1+\binom{n}{1}+(5-\dim_\k M_2)\binom{n}{2}.
		\]
		\item[(2)] If $\vp$ is an automorphism, then $A$ is associated either with a $2$-admissible singleton $\M = (M_2)$, in which case
		\[
		c_n^\vp(A)= 1+\binom{n}{1}+(5-\dim_\k M_2)\binom{n}{2},
		\]
		or with the unique $2$-admissible pair $\P = (\Ga_1^\vp, V_{(1),(1)} \oplus V_{\es,(2)} \oplus V_{\es,(1^2)})$, in which case
		\[
		 c_n^\vp(A) = 1 + \binom{n}{2}.
		\]
	\end{enumerate}
\end{prop}
	
\subsection{Cubic graded codimension sequences of unital superalgebras}
	In this subsection, $\vp$ denotes an automorphism and all $\vp$-algebras are superalgebras. To distinguish this case from the antiautomorphism setting, we fix the following notation:
	\[
	\T^\vp = \T_2, \quad \Id^\vp = \Id^{gr}, \quad c^\vp = c^{gr}, \quad \Ga^\vp = \Ga^{gr}.
	\]	
	To classify $\T_2$-ideals with cubic graded codimension growth, it suffices to classify all $3$-admissible sequences.
	\begin{rmk}\label{rmk2.6}
		Let $\M=(M_3)$ be a $3$-admissible singleton. Then $M_3$ is a proper submodule of $\Ga_3^{gr}$. Since $\Ga_3^{gr}$ has infinitely many submodules, there exist infinitely many $3$-admissible singletons. However, the corresponding isomorphism classes of $M_3$ as $\k\H_3$-modules are finite. Therefore, in the superalgebra setting, our classification of $3$-admissible sequences is carried out up to componentwise isomorphism of $\k\H_3$-modules, which is sufficient for determining the possible graded codimension sequences.
	\end{rmk}
	
	We now turn to the $3$-admissible pairs $\P = (P_2, P_3)$.
	
	The free algebra $\k \lr{Y,Z}$ carries a natural superalgebra structure: let $\K^+$ and $\K^-$ denote the subspaces spanned by all monomials with an even and odd number of $Z$-variables, respectively. Then $\k\lr{Y,Z} = \K^+ \oplus \K^-$. For a polynomial $f(y_1, \dots, y_n, z_1, \dots, z_m)$, denote by $\lr{f}_{\T_2}$ the $\T_2$-ideal of the free superalgebra $\k\lr{Y,Z}$ generated by the polynomial $f$. Every element of $\lr{f}_{\T_2}$ is a linear combination of elements of the form
	\[
	h_l f(g_1^+, \dots, g_n^+, g_1^-, \dots, g_m^-) h_r,
	\]
	where $g_1^+, \dots, g_n^+ \in \K^+$, $g_1^-, \dots, g_m^- \in \K^-$, and $h_l, h_r \in \k\lr{Y,Z}$.
	
\begin{lem}\label{lem2.7}
	Let $P_2$ be an irreducible submodule of $\Ga_2^{gr} = V_{(1^2), \es} \oplus  V_{(1),(1)} \oplus V_{\es,(2)} \oplus V_{\es,(1^2)}$.
	\begin{enumerate}
		\item[(1)]
		If $P_2=V_{(1^2), \es}$, then $\chi({\lr{V_{(1^2), \es}}_{\T_2} \cap \Ga_3^{gr}}) = \chi_{(2,1),\es} + 2\chi_{(1^2),(1)} + \chi_{(1),(2)} + \chi_{(1),(1^2)}$.
		\item[(2)]
		If $P_2=V_{(1), (1)}$, then $\chi({\lr{V_{(1), (1)}}_{\T_2} \cap \Ga_3^{gr}}) = \chi_{(2),(1)} + 2\chi_{(1^2),(1)} + 2\chi_{(1),(2)} + 2\chi_{(1),(1^2)} + 2\chi_{\es, (2,1)}$.
		\item[(3)]
		If $P_2=V_{\es, (2)}$, then $\chi({\lr{V_{\es, (2)}}_{\T_2} \cap \Ga_3^{gr}}) =2\chi_{(1),(2)} + 2\chi_{(1),(1^2)} + \chi_{\es, (3)} + 2\chi_{\es, (2,1)}$.
		\item[(4)]
		If $P_2=V_{\es, (1^2)}$, then $\chi({\lr{V_{\es, (1^2)}}_{\T_2} \cap \Ga_3^{gr}}) =2\chi_{(1),(2)} + 2\chi_{(1),(1^2)} + \chi_{\es, (1^3)} + 2\chi_{\es, (2,1)}$.
	\end{enumerate}
\end{lem}	
	
\begin{proof}
	To compute the $\T_2$-ideal generated by $P_2$, we use the following bases of the homogeneous components of the free superalgebra:
	\[
	\K^+ \cap V_2^{gr} = \sp_\k \{y_1y_2, y_2y_1, z_1z_2, z_2z_1\}, \quad \K^- \cap V_2^{gr} = \sp_\k \{y_1z_2, y_2z_1, z_1y_2, z_2y_1\}.
	\]
	
	(1) Assume $P_2=V_{(1^2), \es}=\k[y_2,y_1]$. Since each symmetric variable $y \in \K^+$, substituting $y$ by the above basis elements of $\K^+ \cap V_2^{gr}$ produces monomials of the form $yy$ and $zz$. Together with the two-sided ideal property, the left $\k\H_3$-module $\lr{V_{(1^2), \es}}_{\T_2} \cap V_3^{gr}$ is generated by elements of the following forms:
	\[
	[yy,y],\quad y[y,y],\quad [y,y]y,\qquad z[y,y],\quad [y,y]z,\qquad [zz,y].
	\]
	We now analyze the $\H_3$-character according to the bipartition type.
	
	\textbf{Case 1.} Bipartition of type $(3,0)$. Only the Specht basis elements $[y_3,y_1,y_2]$, $[y_3,y_2,y_1]$ of $\Ga_3^{gr}$ are involved. Since $[y_3,y_1,y_2]=[y_3,y_1]y_2-y_2[y_3,y_1]$, the commutator $[y,y,y]$ can be expressed in terms of $[y,y]y$ and $y[y,y]$. Hence, the $\H_3$-character $\chi(\lr{V_{(1^2), \es}}_{\T_2} \cap \Ga_3^{gr})$ contains one copy of $\chi_{(2,1),\es}$.
	
	\textbf{Case 2.} Bipartition of type $(2,1)$. By the representation theory of the symmetric group, it suffices to consider the $\S_{\{1,2\}} \times \S_{\{3\}}$-characters. The relevant Specht basis elements are $[z_3,y_1,y_2]$, $[z_3,y_2,y_1]$, $z_3[y_2,y_1]$. By linear independence, the generator $z[y,y]$ yields only $z_3[y_2,y_1]$,  while $[y,y]z$ yields only $[y_2,y_1]z_3$. Note that $[y_2,y_1]z_3=([z_3,y_1,y_2]-[z_3,y_2,y_1])+z_3[y_2,y_1]$. It follows that 
	\[
	z_3[y_2,y_1],\quad [y_2,y_1]z_3 \in \Ga_3^{gr}.
	\]
	The transposition $(12)$ acts by
	\[
	(12)\cdot z_3[y_2,y_1]=-z_3[y_2,y_1], \quad (12)\cdot [y_2,y_1]z_3=-[y_2,y_1]z_3.
	\]
	Thus the $\H_3$-character $\chi(\lr{V_{(1^2), \es}}_{\T_2} \cap \Ga_3^{gr})$ contains two copies of $\chi_{(1^2),(1)}$.
	
	\textbf{Case 3.} Bipartition of type $(1,2)$. It suffices to consider the $\S_{\{1\}} \times \S_{\{2,3\}}$-characters. The relevant Specht basis elements are $[z_3,y_1,z_2]$, $[z_3,z_2,y_1]$, $z_2[z_3,y_1]$ and $z_3[z_2,y_1]$. The generator $[zz,y]$ yields $[z_2z_3,y_1]$ and $[z_3z_2,y_1]$. Note that
	\begin{equation*}
		\begin{split}
		[z_2z_3,y_1]&=[z_3,y_1,z_2]-[z_3,z_2,y_1]+z_2[z_3,y_1]+z_3[z_2,y_1] \in \Ga_3^{gr}, \\
		[z_3z_2,y_1]&=[z_3,y_1,z_2]+z_2[z_3,y_1]+z_3[z_2,y_1] \in \Ga_3^{gr}.
		\end{split}
	\end{equation*}
	The transposition $(23)$ acts by
	\begin{equation*}
		\begin{split}
		(23)\cdot([z_2z_3,y_1]+[z_3z_2,y_1])&=[z_2z_3,y_1]+[z_3z_2,y_1], \\
		(23)\cdot([z_2z_3,y_1]-[z_3z_2,y_1])&=-([z_2z_3,y_1]-[z_3z_2,y_1]).
		\end{split}
	\end{equation*}
	The $\H_3$-character $\chi(\lr{V_{(1^2), \es}}_{\T_2} \cap \Ga_3^{gr})$ contains one copy of $\chi_{(1),(2)}$ and one copy of $\chi_{(1),(1^2)}$.
	
	Combining the three cases we obtain  $\chi({\lr{V_{(1^2), \es}}_{\T_2} \cap \Ga_3^{gr}}) = \chi_{(2,1),\es} + 2\chi_{(1^2),(1)} + \chi_{(1),(2)} + \chi_{(1),(1^2)}$.

	The proofs for the remaining three cases follow by completely analogous character calculations for the generated ideals.
\end{proof}	
	
	Let $\P=(P_2, P_3)$ be a $3$-admissible pair. By Lemma~\ref{lem2.7}, we compute the $\H_3$-character $\chi\bigl( \langle P_2 \rangle_{\T_2} \cap \Ga_3^{gr} \bigr)$, and denote by $\#\chi(P_3)$ the number of possible $\H_3$-characters of $P_3$ for a fixed $P_2$.
	
\begin{table}[H] 
	\caption{Classification of $3$-admissible pairs for unital superalgebras}
	\label{tab1}
	\small
	\centering
	\begin{tabular}{cccccccccc}
		\hline
		$P_2$   & \multicolumn{8}{c}{Irreducible multiplicities in $\chi(\lr{P_2}_{\T_2} \cap \Ga_3^{gr})$} & $\#\chi(P_3)$ \\
		&$(2,1),\es$	&$(2),(1)$	&$(1^2),(1)$   &$(1),(2)$   &$(1),(1^2)$   &$\es,(3)$  &$\es,(1^3)$   &$\es,(2,1)$ &  \\ 
		\hline 
		$V_{(1^2),\es}$   					&1 &0 &2 &1 &1 &0 &0 &0  & 95 \\ 
		$V_{\es,(2)}$             			&0 &0 &0 &2 &2 &1 &0 &2  & 23 \\ 
		$V_{\es,(1^2)}$  						&0 &0 &0 &2 &2 &0 &1 &2  & 23 \\ 
		$V_{(1),(1)}$    							&0 &1 &2 &2 &2 &0 &0 &2  & 7 \\ 
		$V_{\es,(2)} \oplus V_{\es,(1^2)}$   	&0 &0 &0 &2 &2 &1 &1 &2  & 11 \\ 		
		$V_{(1^2),\es} \oplus V_{(1),(1)}$     		&1 &1 &2 &2 &2 &0 &0 &2  & 3 \\ 
		$V_{(1^2),\es} \oplus V_{\es,(2)}$   	&1 &0 &2 &2 &2 &1 &0 &2  & 3 \\ 
		$V_{(1^2),\es} \oplus V_{\es,(1^2)}$ 	&1 &0 &2 &2 &2 &0 &1 &2  & 3 \\ 
		$V_{(1),(1)} \oplus V_{\es,(2)}$   			&0 &1 &2 &2 &2 &1 &0 &2  & 3 \\ 
		$V_{(1),(1)} \oplus V_{\es,(1^2)} $   		&0 &1 &2 &2 &2 &0 &1 &2  & 3 \\ 
		$V_{(1^2),\es} \oplus V_{(1),(1)} \oplus V_{\es,(2)}$   		&1 &1 &2 &2 &2 &1 &0 &2  & 1 \\ 
		$V_{(1^2),\es} \oplus V_{(1),(1)} \oplus V_{\es,(1^2)}$   	&1 &1 &2 &2 &2 &0 &1 &2  & 1 \\ 
		$V_{(1^2),\es} \oplus V_{\es,(2)} \oplus V_{\es,(1^2)}$ &1 &0 &2 &2 &2 &1 &1 &2  & 1 \\ 
		$V_{(1),(1)} \oplus V_{\es,(2)} \oplus V_{\es,(1^2)}$   		&0 &1 &2 &2 &2 &1 &1 &2  & 1 \\ 
		\hline		
	\end{tabular}
\end{table}	
	
\begin{rmk}\label{rmk2.8}
	 Table~\ref{tab1} provides a lower bound on the $\H_3$-character of $P_3$ for any $3$-admissible pair $\P = (P_2, P_3)$. For instance, if $P_2 = V_{(1^2),\es}$, then
	\[
	\chi( \lr{V_{(1^2), \es}}_{\T_2} \cap \Ga_3^{gr} ) \leftrightarrow (1, 0, 2, 1, 1, 0, 0, 0), \quad \chi(\Ga_3^{gr}) \leftrightarrow (1, 1, 2, 2, 2, 1, 1, 2),
	\]
	where each coordinate of these sequences records the multiplicity of the corresponding irreducible $\H_3$-characters. Componentwise subtraction yields the sequence $(0, 1, 0, 1, 1, 1, 1, 2)$. Since group algebra $\k\H_3$ is semisimple, submodules of $\Ga_3^{gr}$ are determined up to isomorphism by the multiplicities of irreducible characters, and these multiplicities can be chosen independently within the given ranges. Hence the number of distinct admissible $\chi(P_3)$ equals
	\[
	1 \cdot 2 \cdot 1 \cdot 2 \cdot 2 \cdot 2 \cdot 2 \cdot 3 - 1 = 95,
	\]
	where the term $-1$ corresponds to the maximal choice $P_3 = \Ga_3^{gr}$. Hence, the corresponding $3$-admissible pair is determined up to componentwise isomorphism of $\k\H_3$-modules.
\end{rmk}

\begin{rmk}\label{rmk2.9}
	Example~\ref{exm2.4} shows that $(\Ga_1^{gr}, V_{(1),(1)} \oplus V_{\es,(2)} \oplus V_{\es,(1^2)})$ is the unique $2$-admissible pair. Moreover, the last row of Table~\ref{tab1} implies that $\lr{V_{(1),(1)} \oplus V_{\es,(2)} \oplus V_{\es,(1^2)}}_{\T_2} \cap \Ga_3^{gr}$ is a maximal submodule of $\Ga_3^{gr}$. Since every polynomial generated by the skew variable $z$ always contains at least a factor $z$, the submodule $V_{(2,1),\es}$, which is spanned by commutators involving only the symmetric variable $y$, satisfies  $V_{(2,1),\es} \not\subset \lr{\Ga_1^{gr}}_{\T_2} \cap \Ga_3^{gr}$. Consequently, the admissible $P_3$ must coincide with the above maximal submodule, and the corresponding $3$-admissible triple is unique. 
\end{rmk}

	The classification of $3$-admissible sequences obtained above yields the complete classification of cubic graded codimension sequences.

\begin{thm}\label{thm2.10}
	Let $A$ be a unital superalgebra whose graded codimension sequence $c_n^{gr}(A)$ has cubic growth, and let $\M$ be the associated $3$-admissible sequence. There are exactly $59$ distinct graded codimension sequences, listed below:
	\begin{enumerate}
		\item[(1)]
		If $\M=(M_3)$ is the $3$-admissible singleton with $M_3$ a proper submodule of $\Ga_3^{gr}$, then:
		\[
		c_n^{gr}(A)= 1+\binom{n}{1}+5\binom{n}{2}+u\binom{n}{3}
		\]
		where $u \in [29]$.
		\item[(2)]
		If $\M=(M_2,M_3)$ is the $3$-admissible pair, then one of the following holds:
		\begin{enumerate}
			\item[(i)]
			\[
			c_n^{gr}(A)= 1+\binom{n}{1}+4\binom{n}{2}+a\binom{n}{3}
			\]
			where $a \in [15]$.
			\item[(ii)]
			\[
			c_n^{gr}(A)= 1+\binom{n}{1}+3\binom{n}{2}+b\binom{n}{3}
			\]
			where $b \in [11] \setminus \{7,10\}$.
			\item[(iii)]
			\[
			c_n^{gr}(A)= 1+\binom{n}{1}+2\binom{n}{2}+c\binom{n}{3}
			\]
			where $c \in [3]$.
			\item[(iv)]
			\[
			c_n^{gr}(A)= 1+\binom{n}{1}+\binom{n}{2}+d\binom{n}{3}
			\]
			where $d \in [2]$.
		\end{enumerate}
		\item[(3)]
		If $\M=(M_1,M_2,M_3)$ is the $3$-admissible triple, then:
		\[
		c_n^{gr}(A)= 1+\binom{n}{2}+2\binom{n}{3}.
		\]
	\end{enumerate}
\end{thm}

\begin{proof}
	All statements follow from Lemma \ref{lem2.3} together with the classification of $3$-admissible sequences. We only prove case (2)(ii), the rest are analogous.
	
	Let $\M=(M_2,M_3)$ be a $3$-admissible pair with $\dim_{\k} M_2 = 2$. Then $M_2$ must be one of the following four submodules:
	\[
	V_{(1),(1)},\quad V_{(1^2),\es} \oplus V_{\es,(2)},\quad V_{(1^2),\es} \oplus V_{\es,(1^2)},\quad V_{\es,(2)} \oplus V_{\es,(1^2)}.
	\]
	
	If $M_2 = V_{(1),(1)}$, by the fourth row of Table \ref{tab1},
	\[
	\chi({\lr{V_{(1),(1)}}_{\T_2} \cap \Ga_3^{gr}}) = \chi_{(2),(1)} + 2\chi_{(1^2),(1)} + 2\chi_{(1),(2)} + 2\chi_{(1),(1^2)} + 2\chi_{\es, (2,1)}.
	\]
	Subtracting this from $\chi(\Ga_3^{gr})$ yields
	\[
	\chi(\Ga_3^{gr})-\chi(\lr{V_{(1),(1)}}_{\T_2} \cap \Ga_3^{gr})=\chi_{(2,1),\es}+\chi_{\es, (3)}+\chi_{\es, (1^3)}.
	\]
	Hence there are exactly $7$ possible $\H_3$-characters for $M_3$, corresponding to dimensions $\dim_{\k} M_3 \in \{25, 26, 27, 28\}$.

	If $M_2=V_{(1^2),\es} \oplus V_{\es,(2)}$, by the seventh  row of Table \ref{tab1}, 
	\[
	\chi({\lr{M_2}_{\T_2} \cap \Ga_3^{gr}}) =\chi_{(2,1),\es} + 2\chi_{(1^2),(1)} + 2\chi_{(1),(2)} + 2\chi_{(1),(1^2)} + \chi_{\es,(3)} + 2\chi_{\es, (2,1)}.
	\]
	Since $\chi(\Ga_3^{gr})-\chi(\lr{M_2}_{\T_2} \cap \Ga_3^{gr})=\chi_{(2),(1)}+\chi_{\es, (1^3)}$, there are exactly $3$ choices for $\chi(M_3)$, with $\dim_{\k} M_3 \in \{25, 26, 28\}$.
	
	If $M_2=V_{(1^2),\es} \oplus V_{\es,(1^2)}$, by the eighth row of Table \ref{tab1},
	\[
	\chi({\lr{M_2}_{\T_2} \cap \Ga_3^{gr}}) =\chi_{(2,1),\es} +  2\chi_{(1^2),(1)} + 2\chi_{(1),(2)} + 2\chi_{(1),(1^2)} + \chi_{\es,(1^3)} + 2\chi_{\es, (2,1)}
	\]
	Since $\chi(\Ga_3^{gr})-\chi(\lr{M_2}_{\T_2} \cap \Ga_3^{gr})=\chi_{(2),(1)}+\chi_{\es, (3)}$, there are again $3$ possibilities for $\chi(M_3)$, with $\dim_{\k}M_3 \in \{25,26,28\}$.
	
	For $M_2=V_{\es,(2)} \oplus V_{\es,(1^2)}$, by the fifth row of Table \ref{tab1}, 
	\[
	\chi({\lr{M_2}_{\T_2} \cap \Ga_3^{gr}}) = 2\chi_{(1),(2)} + 2\chi_{(1),(1^2)} + \chi_{\es,(3)} + \chi_{\es,(1^3)} + 2\chi_{\es, (2,1)}
	\]
	Since $\chi(\Ga_3^{gr})-\chi(\lr{M_2}_{\T_2} \cap \Ga_3^{gr})=\chi_{(2,1),\es}+\chi_{(2),(1)}+2\chi_{(1^2), (1)}$, there are exactly $11$ admissible choices for $\chi(M_3)$, with $\dim_{\k}M_3 \in \{18,20,21,23,24,26,27\}$.
	
	Taking the four cases together, when $\dim_{\k} M_2 = 2$ the possible dimensions of $M_3$ are exactly 
	\[
	\{18, 20, 21, 23, 24, 25, 26, 27, 28\},
	\]
	and the resulting graded codimension sequence takes the form
	\[
	c_n^{gr}(A)= 1+\binom{n}{1}+3\binom{n}{2}+b\binom{n}{3}
	\]
	where $b \in [11] \setminus \{7,10\}$.
\end{proof}

\subsection{Cubic $*$-codimension sequences of unital $*$-algebras}
	In this subsection, $\vp$ denotes an antiautomorphism and all $\vp$-algebras are $*$-algebras. To distinguish this case from the automorphism setting, we fix the following notation:
	\[
	\T^\vp = \T^*, \quad \Id^\vp = \Id^*, \quad c^\vp = c^*, \quad \Ga^\vp = \Ga^*.
	\]

	As in the superalgebra case, the classification of $\T^*$-ideals with cubic $*$-codimension sequences reduces to the classification of $3$-admissible sequences. However, an additional restriction arises: each admissible submodule $P_i$ must be closed under the involution $*$, i.e., $*(P_i) \subset P_i$. 
	
	A $3$-admissible singleton is a sequence $\M = (M_3)$ where $M_3$ is a proper $\k\H_3$-submodule of $\Ga_3^*$ closed under the involution $*$, i.e., $*(M_3) \subset M_3 $. Since the involution $*$ is a linear map, it suffices to verify $*$-invariance on each irreducible summand. Recall the $\H_3$-character of $\Ga_3^*$ is
	\[
	\chi(\Ga_3^*) = \chi_{(2,1),\es} + \chi_{(2),(1)} + 2\chi_{(1^2),(1)} + 2\chi_{(1),(2)} + 2\chi_{(1),(1^2)} + \chi_{\es, (3)} + 2\chi_{\es, (2,1)} + \chi_{\es, (1^3)}.
	\]  
	For an irreducible submodule that occurs with multiplicity one, the $*$-invariance is easily verified. For the characters of multiplicity two, the eigenvectors of the involution $*$ corresponding to the eigenvalue $-1$ are listed below:
	\begin{align*}
		f_{(1^2),(1)} &= [z_3,y_1,y_2]-[z_3,y_2,y_1],\\
		f_{(1),(2)} &= (2[z_3,y_1,z_2]-[z_3,z_2,y_1])+2(z_2[z_3,y_1]+z_3[z_2,y_1]),\\
		f_{(1),(1^2)} &= [z_3,z_2,y_1]+2(z_2[z_3,y_1]-z_3[z_2,y_1]),\\
		f_{\es, (2,1)}^1 &= [z_3,z_1,z_2],\\
		f_{\es, (2,1)}^2 &= [z_3,z_2,z_1];
	\end{align*}
	the eigenvectors of the involution $*$ corresponding to the eigenvalue $1$:
	\begin{align*}
		g_{(1^2),(1)} &= [z_3,y_1,y_2]-[z_3,y_2,y_1]+2z_3[y_2,y_1],\\
		g_{(1),(2)} &= 2[z_3,y_1,z_2]-[z_3,z_2,y_1],\\
		g_{(1),(1^2)} &= [z_3,z_2,y_1],\\
		g_{\es, (2,1)}^1 &= 3[z_3,z_1,z_2]+2(z_1[z_3,z_2]+2z_2[z_3,z_1]+z_3[z_2,z_1]),\\
		g_{\es, (2,1)}^2 &= 3[z_3,z_2,z_1]+2(2z_1[z_3,z_2]+z_2[z_3,z_1]-z_3[z_2,z_1]).
	\end{align*}
	Denote by $U_{\la,\mu}^f$ and $U_{\la,\mu}^g$ the $\k\H_3$-submodules of $\Ga_3^*$ generated by $f_{\la,\mu}$ and $g_{\la,\mu}$, respectively. These are irreducible and correspond to the character $\chi_{\la,\mu}$. Let $M_{\la,\mu}$ be an irreducible $\k\H_3$-submodule of $\Ga_3^*$ with character $\chi_{\la,\mu}$. If $M_{\la,\mu}$ is different from both $U_{\la,\mu}^f$ and $U_{\la,\mu}^g$, then $\langle M_{\la,\mu}\rangle_{\T^*} \cap \Ga_3^* = U_{\la,\mu}^f \oplus U_{\la,\mu}^g$. 
	
\begin{rmk}\label{rmk2.11}
	The character $\chi(\Ga_3^*)$ has four irreducible constituents of multiplicity one and four of multiplicity two. For a multiplicity-one constituent, the corresponding isotypic component is irreducible, so its only $*$-invariant submodules are $0$ and the component itself. For a multiplicity-two constituent, the eigenspace decomposition of $*$ on the isotypic component yields exactly four $*$-invariant submodules: $0$, the $+1$-eigenspace, the $-1$-eigenspace, and the whole isotypic component. Hence the total number of $3$-admissible singletons, i.e., proper $*$-invariant submodules of $\Ga_3^*$, equals
	\[
	2^4\cdot 4^4 - 1 = 4095.
	\]
	Thus, in contrast to the superalgebra case, which is classified up to componentwise isomorphism of $\k\H_3$-modules, this classification is a complete enumeration of concrete submodules.
\end{rmk}

	We now turn to classifying $3$-admissible pairs $\P = (P_2, P_3)$.
	
\begin{lem}\label{lem2.12}
	Let $P_2$ be an irreducible submodule of $\Ga_2^*= V_{(1^2), \es} \oplus  V_{(1),(1)} \oplus V_{\es,(2)} \oplus V_{\es,(1^2)}.$
	\begin{enumerate}
		\item[(1)]
		If $P_2=V_{(1^2), \es}$, then $\chi({\lr{V_{(1^2), \es}}_{\T^*} \cap \Ga_3^*}) = \chi_{(2,1),\es} +\chi_{(2),(1)} + 2\chi_{(1^2),(1)} + \chi_{(1),(2)}$.
		\item[(2)]
		If $P_2=V_{(1), (1)}$, then $\chi({\lr{V_{(1), (1)}}_{\T^*} \cap \Ga_3^*}) = \chi_{(2,1),\es} + \chi_{(2),(1)} + 2\chi_{(1^2),(1)} + 2\chi_{(1),(2)} + 2\chi_{(1),(1^2)} + \chi_{\es, (2,1)}$.
		\item[(3)]
		If $P_2=V_{\es, (2)}$, then $\chi({\lr{V_{\es, (2)}}_{\T^*} \cap \Ga_3^*}) = \chi_{(1^2),(1)} + 2\chi_{(1),(2)} + \chi_{(1),(1^2)} + \chi_{\es, (3)} + \chi_{\es, (1^3)} + 2\chi_{\es, (2,1)}$.
		\item[(4)]
		If $P_2=V_{\es, (1^2)}$, then $\chi({\lr{V_{\es, (1^2)}}_{\T^*} \cap \Ga_3^*}) = \chi_{(1^2),(1)} + \chi_{(1),(2)} + 2\chi_{(1),(1^2)} + \chi_{\es, (1^3)} + 2\chi_{\es, (2,1)}$.
	\end{enumerate}
\end{lem}	
	
\begin{proof}
	Recall that bases for $\K^+ \cap V_2^*$ and $\K^- \cap V_2^*$ are given respectively by 
	\[
	\{y_1 \circ y_2, z_1 \circ z_2, [y_1,z_2], [z_1,y_2]\} \quad \text{and} \quad \{[y_1,y_2], [z_1,z_2], y_1 \circ z_2, z_1 \circ y_2\}.
	\]
	
	(1) Assume $P_2 = V_{(1^2), \es} = \sp_\k \{[y_2, y_1] \}$. Using the above bases of $\K^+ \cap V_2^*$ and $\K^- \cap V_2^*$, the left $\k\H_3$-module $\lr{V_{(1^2), \es}}_{\T^*} \cap V_3^*$ is generated by elements of the following forms:
	\[
	[y \circ y, y], \quad y[y,y], \quad [y,y]y, \qquad [[z,y],y], \quad z[y,y], \quad [y,y]z, \qquad [z \circ z, y].
	\]
	
	\textbf{Case 1.} Bipartition of type $(3,0)$. Combining the identity $[y_3,y_1,y_2]=[y_3,y_1]y_2-y_2[y_3,y_1]$ with the relevant Specht basis elements $[y_3,y_1,y_2]$ and $[y_3,y_2,y_1]$ of $\Ga_3^*$, we see that the $\H_3$-character $\chi\big(\lr{V_{(1^2), \es}}_{\T^*} \cap \Ga_3^*\big)$ contains exactly one copy of $\chi_{(2,1),\es}$.
	
	\textbf{Case 2.} Bipartition of type $(2,1)$. By the representation theory of symmetric groups it suffices to consider the $\S_2\times\S_1$-action. The relevant Specht basis elements of $\Ga_3^*$ are $[z_3,y_1,y_2]$, $[z_3,y_2,y_1]$ and $z_3[y_2,y_1]$, all of which are readily obtained from generators $[z,y,y]$ and $z[y,y]$. The transposition $(12)$ acts by
	\begin{align*}
		(12)\cdot([z_3,y_1,y_2]+[z_3,y_2,y_1])&=([z_3,y_1,y_2]+[z_3,y_2,y_1]),\\
		(12)\cdot([z_3,y_1,y_2]-[z_3,y_2,y_1])&=-([z_3,y_1,y_2]-[z_3,y_2,y_1]),\\
		(12)\cdot z_3[y_2,y_1]&=-z_3[y_2,y_1].
	\end{align*}
	Hence, the $\H_3$-character $\chi( \lr{V_{(1^2), \es}}_{\T^*} \cap \Ga_3^* )$ contains one copy of $\chi_{(2),(1)}$ and two copies of $\chi_{(1^2),(1)}$.	
	
	\textbf{Case 3.} Bipartition of type $(1,2)$. Here it is enough to consider the $\S_1\times\S_2$-action. The relevant Specht basis elements of $\Ga_3^*$ are $[z_3,y_1,z_2]$, $[z_3,z_2,y_1]$, $z_2[z_3,y_1]$ and $z_3[z_2,y_1]$. The generator $[z\circ z,y]$ produces only $[z_2\circ z_3,y_1]$. Note that
	\[
	[z_2 \circ z_3,y_1]=2[z_3,y_1,z_2]-[z_3,z_2,y_1]+2z_2[z_3,y_1]+2z_3[z_2,y_1] = f_{(1),(2)} \in \Ga_3^*,
	\]
	Since $(23)\cdot [z_2 \circ z_3,y_1] = [z_2 \circ z_3,y_1]$, the $\H_3$-character $\chi({\lr{V_{(1^2), \es}}_{\T^*} \cap \Ga_3^*})$ contains one copy of $\chi_{(1),(2)}$.
	
	Combining the three cases we obtain
	\[
	\chi({\lr{V_{(1^2), \es}}_{\T^*} \cap \Ga_3^*}) = \chi_{(2,1),\es} +\chi_{(2),(1)} + 2\chi_{(1^2),(1)} + \chi_{(1),(2)}.
	\] 
	
	In parts (2), (3) and (4), we only detail the bipartitions that appear with multiplicity $2$ in $\chi(\Ga_3^*)$ but with multiplicity $1$ in the corresponding $\chi(\langle P_2\rangle_{\T^*} \cap \Ga_3^*)$. The remaining bipartitions follow by the same arguments as in (1).
	
	(2) Assume $P_2 = V_{(1), (1)} = \sp_\k \{[z_2,y_1], [y_2,z_1] \}$. The $\k\H_3$-module $\lr{V_{(1),(1)}}_{\T^*} \cap V_3^*$ is generated by
	\[
	[[y,y],y],\enspace\enspace [z\circ y,y],\enspace [z,y\circ y],\enspace y[z,y],\enspace [z,y]y,\enspace\enspace [[z,z],y],\enspace [z,[z,y]],\enspace  z[z,y],\enspace  [z,y]z,\enspace\enspace [z,z\circ z].
	\] 
	We concentrate on the bipartition $(\es, (2,1))$, for which the only relevant generator is $[z,z\circ z]$. Observe that $[z_1,z_2\circ z_3]$ and $(12)\cdot [z_1,z_2\circ z_3]$ are linearly independent. Under the involution $*$,
	\begin{align*}
	*([z_1,z_2\circ z_3])=[z_1,z_2\circ z_3]=-\frac{2}{3}g_{\es,(2,1)}^1+\frac{1}{3}g_{\es,(2,1)}^2,\\
	*((12)\cdot [z_1,z_2\circ z_3])=(12)\cdot [z_1,z_2\circ z_3]=\frac{1}{3}g_{\es,(2,1)}^1-\frac{2}{3}g_{\es,(2,1)}^2.
	\end{align*}
	Consequently, the $\H_3$-character $\chi(\lr{V_{(1),(1)}}_{\T^*} \cap \Ga_3^*)$ contains one copy of $\chi_{\es,(2,1)}$.
	
	(3) Assume $P_2 = V_{\es, (2)} = \sp_\k \{z_1 \circ z_2 \}$. The $\k\H_3$-module $\lr{V_{\es, (2)}}_{\T^*} \cap V_3^*$ is generated by 
	\[
	[y,y]\circ z,\qquad (y\circ z)\circ z,\quad y(z\circ z),\quad (z\circ z)y,\qquad [z,z]\circ z,\quad z(z\circ z),\quad (z\circ z)z.
	\] 
	We first consider the bipartition $((1^2), (1))$. The only relevant element is $[y_1,y_2] \circ z_3$. Since
	\[
	[y_1,y_2]\circ z_3=-([z_3,y_1,y_2]-[z_3,y_2,y_1])-2z_3[y_2,y_1]= -g_{(1^2),(1)} \in \Ga_3^*,
	\]
	the $\H_3$-character $\chi(\lr{V_{\es, (2)}}_{\T^*} \cap \Ga_3^*)$ contains one copy of $\chi_{(1^2),(1)}$.
	
	We now turn to the bipartition $((1), (1^2))$. The three forms $(y\circ z)\circ z$, $y(z\circ z)$ and $(z\circ z)y$ yield the following maximal linearly independent set:
	\[
	(y_1\circ z_2)\circ z_3, \quad (y_1\circ z_3)\circ z_2, \quad y_1(z_2\circ z_3), \quad (z_2\circ z_3)y_1.
	\]
	Examining the action of the transposition (23) shows that
	\[
	(y_1\circ z_2)\circ z_3-(y_1\circ z_3)\circ z_2=[z_3,z_2,y_1] = g_{(1),(1^2)} \in \Ga_3^*
	\]
	Thus $\chi(\lr{V_{\es, (2)}}_{\T^*} \cap \Ga_3^*)$ also  contains one copy of $\chi_{(1),(1^2)}$.
	
	(4) Assume $P_2 = V_{\es, (1^2)} = \sp_\k \{ [z_1,z_2] \}$. The $\k\H_3$-module $\lr{V_{\es, (1^2)}}_{\T^*} \cap V_3^*$ is generated by
	\[
	[[y,y],z],\qquad [y\circ z,z],\quad y[z,z],\quad [z,z]y,\qquad [[z,z],z],\quad z[z,z],\quad [z,z]z
	\]
	We first consider the bipartition $((1^2), (1))$. The only relevant element is $[[y_1,y_2],z_3]$. Since
	\[
	[y_1,y_2,z_3]=[z_3,y_2,y_1]-[z_3,y_1,y_2] = -f_{(1^2),(1)} \in \Ga_3^*
	\]
	the $\H_3$-character $\chi(\lr{V_{\es, (1^2)}}_{\T^*} \cap \Ga_3^*)$ contains one copy of $\chi_{(1^2),(1)}$.
	
	We now turn to the bipartition $((1), (2))$. The three forms $[y\circ z,z]$, $y[z,z]$ and $[z,z]y$ yield the following maximal linearly independent set:
	\[
	[y_1\circ z_2,z_3], \quad [y_1\circ z_3,z_2], \quad y_1[z_3,z_2], \quad [z_3,z_2]y_1.
	\]
	Examining the action of the transposition (23) shows that
	\[
	[y_1\circ z_2,z_3]+[y_1\circ z_3,z_2]=-(2[z_3,y_1,z_2]-[z_3,z_2,y_1])-2(z_2[z_3,y_1]+z_3[z_2,y_1]) = -f_{(1),(2)} \in \Ga_3^*,
	\]
	Thus $\chi(\lr{V_{\es, (1^2)}}_{\T^*} \cap \Ga_3^*)$ also  contains one copy of $\chi_{(1),(2)}$.
\end{proof}	
	
	Let $\P=(P_2, P_3)$ be a $3$-admissible pair. By Lemma~\ref{lem2.12}, we compute the $\H_3$-character $\chi\bigl( \langle P_2 \rangle_{\T^*} \cap \Ga_3^* \bigr)$. For a given $P_2$, let $N(P_2)$ denote the number of submodules $P_3 \subset \Ga_3^*$ for which $(P_2, P_3)$ is a $3$-admissible pair. In contrast to the superalgebra setting, the $*$-invariance condition makes this a count of specific submodules rather than of isomorphism classes.
	
\begin{table}[H]
	\caption{Classification of $3$-admissible pairs for unital $*$-algebras}
	\label{tab2}
	\small
	\centering
	\begin{tabular}{cccccccccc}
		\hline
		$P_2$ & \multicolumn{8}{c}{Irreducible multiplicities in $\chi(\lr{P_2}_{\T^*} \cap \Ga_3^*)$} & $N(P_2)$ \\
		&$(2,1),\es$	&$(2),(1)$	&$(1^2),(1)$   &$(1),(2)$   &$(1),(1^2)$   &$\es,(3)$  &$\es,(1^3)$   &$\es,(2,1)$ &\\  
		\hline 
		$V_{(1^2),\es}$   					&1 &1 &2 &1 &0 &0 &0 &0  & 127 \\ 
		$V_{\es,(1^2)}$  						&0 &0 &1 &1 &2 &0 &1 &2  & 31 \\ 		 
		$V_{\es,(2)}$             			&0 &0 &1 &2 &1 &1 &1 &2  & 15 \\ 
		$V_{(1),(1)}$    							&1 &1 &2 &2 &2 &0 &0 &1  & 7\\ 
		$V_{(1^2),\es} \oplus V_{(1),(1)}$     		&1 &1 &2 &2 &2 &0 &0 &1  & 7 \\ 
		$V_{(1^2),\es} \oplus V_{\es,(1^2)}$ 	&1 &1 &2 &1 &2 &0 &1 &2  & 3 \\ 
		$V_{\es,(2)} \oplus V_{\es,(1^2)}$   	&0 &0 &2 &2 &2 &1 &1 &2  & 3 \\ 
		$V_{(1^2),\es} \oplus V_{\es,(2)}$   	&1 &1 &2 &2 &1 &1 &1 &2 & 1 \\ 
		$V_{(1),(1)} \oplus V_{\es,(1^2)} $   		&1 &1 &2 &2 &2 &0 &1 &2  & 1 \\ 
		$V_{(1^2),\es} \oplus V_{(1),(1)} \oplus V_{\es,(1^2)}$   	&1 &1 &2 &2 &2 &0 &1 &2  & 1 \\
		\hline
		$V_{(1),(1)} \oplus V_{\es,(2)}$   			&1 &1 &2 &2 &2 &1 &1 &2  & 0 \\ 
		$V_{(1^2),\es} \oplus V_{(1),(1)} \oplus V_{\es,(2)}$   		&1 &1 &2 &2 &2 &1 &1 &2  & 0 \\ 
		$V_{(1^2),\es} \oplus V_{\es,(2)} \oplus V_{\es,(1^2)}$ &1 &1 &2 &2 &2 &1 &1 &2  & 0 \\ 
		$V_{(1),(1)} \oplus V_{\es,(2)} \oplus V_{\es,(1^2)}$   		&1 &1 &2 &2 &2 &1 &1 &2  & 0 \\ 
		\hline		
	\end{tabular}
\end{table}
	
	From Table~\ref{tab2}, there are precisely $10$ choices of $P_2$ for which a $3$-admissible pair $(P_2, P_3)$ exists, yielding a total of $196$ such pairs. The last four rows of Table~\ref{tab2} show that for those $P_2$, the generated ideal already saturates $\Ga_3^*$, so no admissible $P_3$ exists. In particular, the final row gives $\langle V_{(1),(1)} \oplus V_{\es,(2)} \oplus V_{\es,(1^2)}\rangle_{\T^*} \cap \Ga_3^* = \Ga_3^*$. Hence no $3$-admissible triple exists in the $*$-algebra case.

	Based on the classification of $3$-admissible sequences, we now state the complete classification of cubic $*$-codimension sequences.

\begin{thm}	\label{thm2.13}
	Let $A$ be a unital $*$-algebra whose $*$-codimension sequence $c_n^*(A)$ has cubic growth, let $\M$ denote its associated $3$-admissible sequence. There are exactly $54$ distinct $*$-codimension sequences, listed below.
	\begin{enumerate}
		\item[(1)]
		If $\M$ is a $3$-admissible singleton, then:
		\[
		c_n^*(A)= 1+\binom{n}{1}+5\binom{n}{2}+u\binom{n}{3}
		\]
		where $u \in [29]$.
		\item[(2)]
		If $\M$ is a $3$-admissible pair, then one of the following holds:
		\begin{enumerate}
			\item[(i)]
			\[
			c_n^*(A)= 1+\binom{n}{1}+4\binom{n}{2}+a\binom{n}{3}
			\]
			where $a \in [15]$.
			\item[(ii)]
			\[
			c_n^*(A)= 1+\binom{n}{1}+3\binom{n}{2}+b\binom{n}{3}
			\]
			where $b \in [5]$. 
			\item[(iii)]
			\[
			c_n^*(A)= 1+\binom{n}{1}+2\binom{n}{2}+c\binom{n}{3}
			\]
			where $c \in [4]$. 
			\item[(iv)]
			\[
			c_n^*(A)= 1+\binom{n}{1}+\binom{n}{2}+\binom{n}{3}.
			\]
		\end{enumerate}
	\end{enumerate}
\end{thm}

\begin{proof}
	The proof is parallel to that of Theorem~\ref{thm2.10}, using Lemma~\ref{lem2.3} together with the classification of admissible singletons in Remark~\ref{rmk2.11} and the data for admissible pairs in Table~\ref{tab2}.
\end{proof}

\section{Generating identities for unital $\vp$-algebras of growth at most quadratic}
\label{sec3}
	In \cite{BFXYZZZ}, the authors investigate minimal-degree single multilinear generators for the $\T$-ideals of unital algebras and refer to such polynomials as generating identities. This section aims to determine the generating identities of unital $\vp$-algebras whose $\vp$-codimension sequences have at most quadratic growth.
		
\begin{Def} \label{def3.1}
	Let $A$ be a $\vp$-algebra and let $\Id^\vp(A)$ be its associated $\T^\vp$-ideal. A multilinear polynomial $f \in V_n^\vp$ is called a \emph{generating identity} of $A$ if $\lr{f}_{\T^{\vp}} = \Id^{\vp}(A)$, and no $\vp$-identity of degree less than $n$ generates $\Id^{\vp}(A)$. The integer $n$ is called the \emph{generating degree} of $A$ and is denoted by $\gd(A)$.
\end{Def}	
	
	The finite basis property for $\vp$-algebras was established in \cite{AB,AGK,Sv11,Sv14}. Moreover, it was proved in \cite{BXYZZ} that every $\T$-ideal of a unital algebra is generated by a single element. Using similar methods, we prove an analogous result for $\vp$-algebras.
	
\begin{lem} \label{lem3.2}
	Every unital $\vp$-algebra over a field of characteristic zero admits a generating identity.
\end{lem}

\begin{proof}
	By the finite basis property, the $\T^{\vp}$-ideal $\Id^\vp(A)$ is generated by finitely many multilinear polynomials. Let $n$ be the maximal degree among these generators. Since $A$ is unital, any generator of degree $k < n$ may be multiplied on the right by the symmetric variables $y_{k+1}, \dots, y_n$, yielding a multilinear polynomial of degree $n$ in $\Id^\vp(A)$. Replacing each auxiliary variable $y_{k+1},\dots,y_n$ by the identity element $1$ recovers the original polynomial. Hence we may assume that all generators have degree exactly $n$.
	
	Under the natural action of $\H_n$, these polynomials span a $\k\H_n$-submodule of $V_n^{\vp}$. The space $V_n^\vp$ itself is a cyclic $\k\H_n$-module, generated for instance by $(y_1+z_1)\cdots(y_n+z_n)$. Since $\k$ has characteristic zero, the group algebra $\k\H_n$ is semisimple, so every submodule of a cyclic module remains cyclic. It follows that $\Id^\vp(A)$ is generated, as a $\T^\vp$-ideal, by a single multilinear polynomial. Hence, every unital $\vp$-algebra over a field of characteristic zero admits a generating identity.
\end{proof}

	The generating degree admits an upper bound derived from \cite{Ma09,MMM}.
	
\begin{lem}\label{lem3.3}\cite{Ma09,MMM}
	Let $k$ and $i$ be positive integers. If $\vp$ is an automorphism and $k$ is odd, then $\Ga_{k+i}^\vp$ is generated by $\Ga_k^\vp$ and $[y_1,y_2]\cdots[y_k,y_{k+1}]$; otherwise, $\Ga_{k+i}^\vp$ is generated by $\Ga_k^\vp$.
\end{lem}

\begin{lem} \label{lem3.4}
	Let $A$ be a $\vp$-algebra and $\M=(M_s,\dots,M_k)$ denote its associated $k$-admissible sequence.
	\begin{enumerate}
		\item[(1)]
		If $\vp$ is an automorphism and $k$ is even, then $s \leq \gd(A) \leq k+2$.
		\item[(2)]
		If either $\vp$ is an automorphism and $k$ is odd, or $\vp$ is an antiautomorphism, then $s \leq \gd(A) \leq k+1$.
		\item[(3)]
		For $i = s,\dots,k$ let $\alpha_i$ be a generator of $M_i$, and for $j = k+1,k+2$ let $\beta_j$ be a generator of $\Ga_j^\vp$.
		    \begin{enumerate}
			\item[(i)]
			If $\gd(A)=w \leq k$, then
			\[
			\Id^\vp(A) = \Bigl\langle \sum_{i=s}^{w} \alpha_i y_{i+1}\cdots y_w \Bigr\rangle_{\T^\vp},
			\]
			where the product $y_{i+1}\cdots y_w$ is understood to be $1$ when $i=w$.
			\item[(ii)]
			If $\gd(A)=k+1$, then
			\[
			\Id^\vp(A) = \Bigl\langle \sum_{i=s}^{k} \alpha_i y_{i+1}\cdots y_{k+1} + \beta_{k+1} \Bigr\rangle_{\T^\vp}.
			\]
			\item[(iii)]
			If $\gd(A)=k+2$, then
			\[
			\Id^\vp(A) = \Bigl\langle \sum_{i=s}^{k} \alpha_i y_{i+1}\cdots y_{k+2} + \beta_{k+1} y_{k+2} + \beta_{k+2} \Bigr\rangle_{\T^\vp}.
			\]	
		\end{enumerate}
	\end{enumerate}
\end{lem}	
	
\begin{proof}
	The admissible sequence $(M_s,\dots,M_k)$ determines the full sequence $(M_s,\dots, M_k, \Ga_{k+1}^\vp, \Ga_{k+2}^\vp, \dots)$. By Lemma~\ref{lem3.3}, all terms beyond $\Ga_{k+2}^\vp$ are generated by $\Ga_{k+1}^\vp$ and $\Ga_{k+2}^\vp$, which yields the bounds on $\gd(A)$ in (1) and (2). The formulas in (3) follow by lifting the lower-degree generators to the generating degree, exactly as in the proof of Lemma~\ref{lem3.2}.
\end{proof}	
	
	Lemma~\ref{lem3.4} yields a general approach to compute generating identities. We next determine the single generators of $\Ga_k^\vp$ as left $\k\H_k$-modules for $k=1,2,3$.
	
\begin{lem}	\label{lem3.5}
	\begin{enumerate}
		\item[(1)]
		If $f_{\Ga_1^\vp}=z_1$, then $\Ga_1^\vp = \k\H_1 \cdot f_{\Ga_1^\vp}$.
		\item[(2)]
		If $f_{\Ga_2^\vp}=[y_2,y_1]+[z_2,y_1]+z_1z_2$, then $\Ga_2^\vp = \k\H_2 \cdot f_{\Ga_2^\vp}$.
		\item[(3)]
		If $f_{\Ga_3^\vp}=[y_3,y_2,y_1]+[y_3,z_2,y_1]+z_3[y_2,y_1]+[y_3,z_2,z_1]+z_3[z_2,y_1]+z_1z_2z_3$, then $\Ga_3^\vp = \k\H_3 \cdot f_{\Ga_3^\vp}$.
	\end{enumerate} 
\end{lem}

\begin{proof}
	Assertion (1) is trivial; we only prove case (2), and assertion (3) follows by parallel reasoning.
	
	It is clear that $f_{\Ga_2^{\vp}} \in \Ga_2^{\vp}$ from the Specht basis for $\Ga_2^{\vp}$. Hence it suffices to verify that $f_{\Ga_2^\vp}$ generates a $5$-dimensional space under the action of $\H_2=\Z_2 \wr \S_2$. The explicit group actions are computed below:
	\[
	\begin{array}{r@{}l}
		(1,1;(1)) \cdot f_{\Ga_2^\vp} =& [y_2,y_1]+[z_2,y_1]+z_1z_2, \\
		(-1,1;(1)) \cdot f_{\Ga_2^\vp} =& [y_2,y_1]+[z_2,y_1]+(-z_1)z_2 = [y_2,y_1]+[z_2,y_1]-z_1z_2, \\
		(1,-1;(1)) \cdot f_{\Ga_2^\vp} =& [y_2,y_1]+[-z_2,y_1]+z_1(-z_2) = [y_2,y_1]-[z_2,y_1]-z_1z_2, \\
		(-1,-1;(1)) \cdot f_{\Ga_2^\vp} =& [y_2,y_1]+[-z_2,y_1]+(-z_1)(-z_2) = [y_2,y_1]-[z_2,y_1]+z_1z_2, \\
		(1,1;(12)) \cdot f_{\Ga_2^\vp} =& [y_1,y_2]+[z_1,y_2]+z_2z_1 = -[y_2,y_1]-[y_2,z_1]+[z_2,z_1]+z_1z_2, \\
		(-1,1;(12)) \cdot f_{\Ga_2^\vp} =& [y_1,y_2]+[-z_1,y_2]+z_2(-z_1) = -[y_2,y_1]+[y_2,z_1]-[z_2,z_1]-z_1z_2, \\
		(1,-1;(12)) \cdot f_{\Ga_2^\vp} =& [y_1,y_2]+[z_1,y_2]+(-z_2)z_1 = -[y_2,y_1]-[y_2,z_1]-[z_2,z_1]-z_1z_2, \\
		(-1,-1;(12)) \cdot f_{\Ga_2^\vp} =& [y_1,y_2]+[-z_1,y_2]+(-z_2)(-z_1) = -[y_2,y_1]+[y_2,z_1]+[z_2,z_1]+z_1z_2. \\
	\end{array}
	\]
	Relative to the Specht basis
	\[
	\big\{[y_2,y_1], [z_2,y_1], [y_2,z_1], [z_2,z_1], z_1z_2\big\}
	\]
	the coefficient matrix of the above vectors is
	\[
	\begin{pmatrix}
		\begin{array}{rrrrrrrr}
			1& 1& 1& 1& -1& -1& -1& -1 \\
			1& 1& -1& -1& 0& 0& 0& 0 \\
			0& 0& 0& 0& -1& 1& -1& 1 \\
			0& 0& 0& 0& 1& -1& -1& 1 \\
			1& -1& -1& 1& 1& -1& -1& 1 \\
		\end{array}
	\end{pmatrix}.
	\]
	A direct computation shows that this matrix has rank $5$.
\end{proof}

\begin{prop} \label{prop3.6}
	Let $A$ be a $\vp$-algebra and let $\M=(M_s,\dots,M_k)$ denote its associated $k$-admissible sequence.
	\begin{enumerate}
		\item[(1)]
	If $k=0$, then the following hold:
		\begin{itemize}
			\item if $\vp$ is an antiautomorphism, then $f_{\Ga_1^\vp} = z_1$ is a generating identity of $A$;
			\item if $\vp$ is an automorphism, then $f_{\Ga_1^\vp}y_2 + [y_2, y_1]$ is a generating identity of $A$.
		\end{itemize}
		\item[(2)]
		If $k = 1$, then $f_{\Ga_2^\vp} = [y_2, y_1] + [z_2, y_1] + z_1z_2$ is a generating identity of $A$.
		\item[(3)]
		If $k = 2$, the corresponding results are collected in the Table~\ref{tab3}. Let $f_4 = f_{\Ga_3^\vp}\cdot y_4 + [y_2,y_1][y_4,y_3]$.
	\end{enumerate}
\end{prop}
	
			
\begin{table}[H]
	\caption{The generating identities with quadratic $\vp$-codimension}
	\label{tab3}
	\small
	\centering
	\begin{tabular}{rrr}
		\hline
		&\multicolumn{2}{c}{Generating identity}\\
		\multicolumn{1}{c}{2-admissible sequence}	&\multicolumn{1}{c}{Automorphism}	&\multicolumn{1}{c}{Antiautomorphism}   \\ 
		\hline 
		$V_{(1^2),\es}$   				   &$[y_2,y_1]y_3+f_{\Ga_3^\vp}$ &$[y_2,y_1]y_3+f_{\Ga_3^\vp}$  \\ 
		$V_{(1),(1)} \oplus V_{(1^2),\es}$  		   &$([z_2,y_1]+[y_2,y_1])y_3+f_{\Ga_3^\vp}$ &$([z_2,y_1]+[y_2,y_1])y_3+f_{\Ga_3^\vp}$   \\ 
		$V_{\es,(2)} \oplus V_{(1^2),\es}$    &$(z_2 \circ z_1+[y_2,y_1])y_3+f_{\Ga_3^\vp}$ &$(z_2 \circ z_1+[y_2,y_1])y_3+f_{\Ga_3^\vp}$   \\ 
		$V_{\es,(1^2)} \oplus V_{(1^2),\es}$ &$([z_2,z_1]+[y_2,y_1])y_3+f_{\Ga_3^\vp}$ &$([z_2,z_1]+[y_2,y_1])y_3+f_{\Ga_3^\vp}$   \\ 
		$V_{(1),(1)} \oplus V_{\es,(1^2)} \oplus V_{(1^2),\es}$ &$([z_2,y_1] + [z_2,z_1] + [y_2,y_1])y_3 + f_{\Ga_3^\vp}$ &\multicolumn{1}{c}{$([z_2,y_1] + [z_2,z_1] + [y_2,y_1])y_3 + f_{\Ga_3^\vp}$}   \\ 
		$V_{(1),(1)} \oplus V_{\es,(2)} \oplus V_{(1^2),\es}$ &$([z_2,y_1] + z_2 \circ z_1 + [y_2,y_1])y_3 + f_{\Ga_3^\vp}$ &\multicolumn{1}{c}{$[z_2,y_1] + z_2 \circ z_1 + [y_2,y_1]$}    \\ 
		$V_{\es,(2)} \oplus V_{\es,(1^2)} \oplus V_{(1^2),\es}$   &$(z_1z_2+[y_2,y_1])y_3 + f_{\Ga_3^\vp}$ &\multicolumn{1}{c}{$z_1z_2+[y_2,y_1]$}   \\ 
		\hline	
		$0$  									&$f_4$ &$f_{\Ga_3^\vp}$   \\ 
		$V_{(1),(1)}$   						&$[z_2,y_1]y_3y_4 + f_4$ &$[z_2,y_1]y_3+f_{\Ga_3^\vp}$   \\ 
		$V_{\es,(2)}$   					&$(z_2 \circ z_1)y_3y_4 + f_4$ &$(z_2 \circ z_1)y_3+f_{\Ga_3^\vp}$  \\ 
		$V_{\es,(1^2)}$   				&$[z_2,z_1]y_3y_4 + f_4$ &$[z_2,z_1]y_3+f_{\Ga_3^\vp}$   \\ 
		$V_{\es,(2)} \oplus V_{\es,(1^2)}$   			&$z_1z_2y_3y_4 + f_4$ &$z_1z_2y_3+f_{\Ga_3^\vp}$  \\ 
		$V_{(1),(1)} \oplus V_{\es,(1^2)}$   	&$([z_2,y_1]+[z_2,z_1])y_3y_4 + f_4$ &$([z_2,y_1]+[z_2,z_1])y_3+f_{\Ga_3^\vp}$  \\ 
		$V_{(1),(1)} \oplus V_{\es,(2)}$   		&$([z_2,y_1]+z_2 \circ z_1)y_3y_4 + f_4$ &\multicolumn{1}{c}{$[z_2,y_1]+z_2 \circ z_1$}   \\ 
		$V_{(1),(1)} \oplus V_{\es,(2)} \oplus V_{\es,(1^2)}$ &$([z_2,y_1]+z_1z_2)y_3y_4 + f_4$ &\multicolumn{1}{c}{$[z_2,y_1]+z_1z_2$}  \\ 
		$(V_{\es,(1)}, V_{(1),(1)} \oplus V_{\es,(2)} \oplus V_{\es,(1^2)})$  &$z_1y_2y_3y_4 + f_4$ &\multicolumn{1}{c}{not applicable}   \\ 
		\hline		
	\end{tabular}
\end{table}

\begin{proof}
	Cases (1) and (2) are elementary: the stated generating identities follow directly from Lemmas~\ref{lem3.4} and \ref{lem3.5}. For (3), Proposition~\ref{prop2.5} classifies all $2$-admissible sequences, which are listed in the first column of Table~\ref{tab3}. In the automorphism case, according to Table~\ref{tab1}, no submodule $M_2$ of $\Ga_2^{gr}$ generates the whole $\Ga_3^{gr}$. Hence, by Lemma~\ref{lem3.3}, the generating degree is $3$ if $M_2$ contains $V_{(1^2),\emptyset}$, and $4$ otherwise. In the antiautomorphism case, the $M_2$ corresponding to the last four rows of Table~\ref{tab2} already generate $\Ga_3^*$, so the generating degree is $2$. For all remaining $M_2$, the generating degree is $3$. The generators of the relevant submodules are obtained by the same method as in Lemma~\ref{lem3.5}. Substituting these into the formulas of Lemma~\ref{lem3.4}(3) yields the generating identities listed in Table~\ref{tab3}.
\end{proof}

\noindent \textbf{Declaration of competing interest}.

	The authors declare that they have no known competing financial interests or personal relationships that could have appeared to influence the work reported in this paper.

\noindent\textbf{Acknowledgments}.

	Yan-Hong Bao was partially supported by the National Natural Science Foundation of China (No. 12371015) and Excellent University Research and Innovation Team in Anhui Province (No. 2024AH010002). Jiang-Nan Xu was supported by the Scientific Research Foundation of Huainan Normal University (No. 826061). Yuan-Feng Zhang was partially supported by the Scientific Research Project of the Education Department of Anhui Province (No. 2025AHGXZK50142).

\makeatletter
\renewcommand{\@biblabel}[1]{\hfill\mbox{[#1]}}


\begin{thebibliography}{9999}
	\raggedbottom
	\bibitem{AB}
	E. Aljadeff, A. Kanel-Belov,
	Representability and Specht problem for  $G$-graded algebras,
	Adv. Math. {\bf 225} (2010), no. 5, 2391–2428.
	
	\bibitem{AGK}
	E. Aljadeff, A. Giambruno, Y. Karasik,
	Polynomial identities with involution, superinvolutions and the Grassmann envelope,
	Proc. Amer. Math. Soc. {\bf 145} (2017), no. 5, 1843–1857.
	
	\bibitem{BFXYZZZ}
	Y.-H. Bao, D.-X. Fu, J.-N. Xu, Y. Ye, J. J. Zhang, Y.-F. Zhang, Z.-B. Zhao,
	Codimension sequence, grade, and generating degree: an operadic approach,
	Available at arXiv: 2504.19168.
	
	\bibitem{BXYZZ}
	Y.-H. Bao, J.-N. Xu, Y. Ye, J. J. Zhang, Y.-F. Zhang,
	Ideals of the associative algebra operad, 
	Available at arXiv: 2307.07669.
	
	\bibitem{BYZ}
	Y.-H. Bao, Y. Ye, J. J. Zhang, 
	Truncation of unitary operads, 
	Adv. Math. {\bf 372} (2020), 107290.
	
	\bibitem{BSSV}
	D. C. L. Bessades, R. B. dos Santos, M. L. O. Santos, A. C. Vieira, 
	Superalgebras and algebras with involution: classifying varieties of quadratic growth,
	Comm. Algebra {\bf 49} (2021), no. 6, 2476–2490.
	
	\bibitem{Dr}
	V. Drensky, 
	Free Algebras and PI-Algebras, Graduate course in algebra,
	Springer-Verlag Singapore, Singapore, 2000. xii+271 pp.
	
	\bibitem{DG}
	V. Drensky, A. Giambruno,
	Cocharacters, codimensions and Hilbert series of the polynomial identities for  $2 \times 2$  matrices with involution,
	Canad. J. Math. {\bf 46} (1994), no. 4, 718–733.

	\bibitem{GM}
	A. Giambruno, D. La Mattina,
	PI-algebras with slow codimension growth,
	J. Algebra {\bf 284} (2005), no. 1, 371–391.
	
	\bibitem{GMM}
	A. Giambruno, D. La Mattina, P. Misso,
	Polynomial identities on superalgebras: classifying linear growth,
	J. Pure Appl. Algebra {\bf 207} (2006), no. 1, 215–240.
	
	\bibitem{GMP}
	A. Giambruno, D. La Mattina, V. M. Petrogradsky, 
	Matrix algebras of polynomial codimension growth, 
	Israel J. Math. {\bf 158} (2007), 367--378.
	
	\bibitem{GMZ}
	A. Giambruno, D. La Mattina,  M. Zaicev, 
	Classifying the minimal varieties of polynomial growth, 
	Canad. J. Math. {\bf 66} (2014), no. 3, 625--640.
	
	\bibitem{GZ}
	A. Giambruno, M. Zaicev,
	Polynomial identities and asymptotic methods.
	Mathematical Surveys and Monographs, 122. American Mathematical Society, Providence, RI, 2005.
	
	\bibitem{GSV}
	T. A. Gouveia, R. B. dos Santos, A. C. Vieira,
	Minimal  $*$-varieties and minimal supervarieties of polynomial growth,
	J. Algebra {\bf 552} (2020), 107–133.
	
	\bibitem{Ke79}
	A. R. Kemer,
	Varieties of finite rank, 
	Proceeding of 15th All the Union Algebraic Conference, Krasnoyarsk, Vol. 2. (in Russian).
	
	\bibitem{Ke87}
	A. R. Kemer,
	Finite basability of identities of associative algebras,
	Algebra i Logika {\bf 26} (1987), no. 5, 597–641, 650.
	
	\bibitem{Ma09}
	D. La Mattina,
	Polynomial codimension growth of graded algebras. 
	In: Groups, Rings and Group Rings. Contemp. Math. vol. 499, Providence, RI, AMS, 2009,  189–197.
	
	\bibitem{MMM}
	D. La Mattina, S. Mauceri, P. Misso,
	Polynomial growth and identities of superalgebras and star-algebras,
	J. Pure Appl. Algebra {\bf 213} (2009), no. 11, 2087–2094.
	
	\bibitem{MM}
	D. La Mattina, P. Misso,
	Algebras with involution with linear codimension growth,
	J. Algebra {\bf 305} (2006), no. 1, 270–291.
	
	\bibitem{OV}
	M. A. de Oliveira, A. C. Vieira,
	Varieties of unitary algebras with small growth of codimensions,
	Internat. J. Algebra Comput. \textbf{31}(2021), no.2, 257--277.
	
	\bibitem{Re}
	A. Regev, 
	Existence of identities in $A\otimes B$,
	Israel J. Math. {\bf 11} (1972), 131--152. 																	
	
	\bibitem{Sp}  
	W. Specht, 
	Gesetze in Ringen. I. (German), 
	Math. Z. {\bf 52} (1950), 557--589.
	
	\bibitem{Sv11}
	I. Sviridova,
	Identities of pi-algebras graded by a finite abelian group,
	Comm. Algebra {\bf 39} (2011), no. 9, 3462–3490.
	
	\bibitem{Sv14}
	I. Sviridova,
	Finite basis problem for identities with involution,
	Available at arXiv: 1410.2233.
	
\end{thebibliography}
\end{document}